\documentclass[12pt]{article}
\usepackage{amsmath,amssymb}
\usepackage{graphicx}
\begin{document}
\newtheorem{theorem}{Theorem}
\newtheorem{lemma}{Lemma}
\newtheorem{corollary}{Corollary}
\newtheorem{proof}{Proof}
\renewcommand{\thetheorem}{\thesection.\arabic{theorem}}
\renewcommand{\thelemma}{\thesection.\arabic{lemma}}
\newcommand{\mysection}[1]{\section{#1}\setcounter{equation}{0}
\setcounter{theorem}{0}} \setcounter{lemma}{0}

\title
{{\bf Generalization of the Secant Method for Nonlinear Equations}\\
(extended version)}
\author
{Avram Sidi\\ Computer Science Department\\
Technion - Israel Institute of Technology\\
Haifa 32000, Israel\\ \\
e-mail: asidi@cs.technion.ac.il\\
http://www.cs.technion.ac.il/$\sim$asidi}
\date{November 2020}

 \maketitle
 \begin{center}
 Original version appeared  in: \ {\em Applied Mathematics E-Notes}, 8:115--123, 2008.
 \end{center}
 \thispagestyle{empty}

\newpage
\begin{abstract}
The  secant method is a very effective numerical procedure used
for solving nonlinear equations of the form $f(x)=0$. It is
derived via a linear interpolation procedure and employs only
values of $f(x)$ at the approximations to the root of $f(x)=0$,
hence it computes $f(x)$ only once per iteration. In this note, we
generalize it by replacing the relevant linear interpolant by a
suitable $(k+1)$-point polynomial of interpolation, where $k$ is
an integer at least 2. Just as the secant method, this
generalization too enjoys the property that it computes $f(x)$
only once per iteration.  We provide its error in closed form and
analyze its order of convergence $s_k$. We show that this order of
convergence is greater than that of the secant method, and it
increases towards $2$ as $k\to \infty$. (Indeed,  $s_7=1.9960\cdots$, for example.)
This is true for the  efficiency index of the method too.
We also confirm the theory via an illustrative example.

\end{abstract}
\thispagestyle{empty}
\newpage
\pagenumbering{arabic}
\section{Introduction} \label{se1}
Let $\alpha$ be the solution to the equation
\begin{equation}\label{eq1} f(x)=0.\end{equation}
 An effective iterative method used for
solving \eqref{eq1} that makes direct use of $f(x)$ [but no
derivatives of $f(x)$] is the {\em secant method} that is
discussed in many books on numerical analysis. See, for example,
Atkinson~\cite{Atkinson:1989:INA},
Henrici~\cite{Henrici:1964:ENA}, Ralston and
Rabinowitz~\cite{Ralston:1978:FCN}, and Stoer and
Bulirsch~\cite{Stoer:2002:INA}. See also the recent note
\cite{Sidi:2006:UTR} by the author, in which the treatment of the
secant method and those of the Newton--Raphson, regula falsi, and
Steffensen methods are presented in a unified manner.

 This method is derived by a linear interpolation
procedure as follows: Starting with two initial approximations
$x_0$ and $x_1$ to the solution $\alpha$ of \eqref{eq1}, we
compute a sequence of approximations $\{x_n\}^\infty_{n=0}$, such
that the approximation $x_{n+1}$ is determined as the point of
intersection (in the $x$-$y$ plane) of the straight line  through
the points $(x_n,f(x_n))$ and $(x_{n-1},f(x_{n-1}))$ with the
$x$-axis.  Since the equation of this straight line is
\begin{equation}\label{eq2}
y=f(x_n)+\frac{f(x_n)-f(x_{n-1})}{x_n-x_{n-1}}\,(x-x_n),
\end{equation} $x_{n+1}$ is given as
\begin{equation}\label{eq3}
x_{n+1}=x_n-\frac{f(x_n)}{\displaystyle\frac{f(x_n)-f(x_{n-1})}{x_n-x_{n-1}}}.\end{equation}
In terms of divided differences, \eqref{eq3} can be written in the
form
\begin{equation}\label{eq4}
x_{n+1}=x_n-\frac{f(x_n)}{f[x_n,x_{n-1}]}.\end{equation}
 Again, in terms of divided differences, the error in $x_{n+1}$
is given as in
\begin{equation}\label{eq5}
x_{n+1}-\alpha=\frac{f[x_n,x_{n-1},\alpha]}{f[x_n,x_{n-1}]}\,(x_n-\alpha)(x_{n-1}-\alpha).
\end{equation}
Then, provided $f(x)$ is twice continuously differentiable in a
closed interval $I$ containing
 $\alpha$ in its interior, and  provided $x_{n-1},x_n\in I$,
 \eqref{eq5} becomes
\begin{gather}
x_{n+1}-\alpha=
\frac{f''(\xi_n)}{2f'(\eta_n)}\,(x_n-\alpha)(x_{n-1}-\alpha),\nonumber\\
\xi_n\in\text{int}(x_n,x_{n-1},\alpha), \quad \eta_n\in
\text{int}(x_n,x_{n-1}). \label{err-secant} \end{gather}
 In case,
$f'(\alpha)\neq 0$ and $x_0$ and $x_1$ are sufficiently close to
$\alpha$, there holds $\lim_{n\to\infty}x_n=~\alpha$, and hence
$$\lim_{n\to\infty}\frac{x_{n+1}-\alpha}{(x_n-\alpha)(x_{n-1}-\alpha)}=
\frac{f''(\alpha)}{2f'(\alpha)}.$$ From this, one derives the
conclusion that the  order of convergence  of the secant method is
at least $(1+\sqrt{5})/2$.

Another way of obtaining the secant method, of interest to us in
the present work,  is via a variation of the Newton--Raphson
method. Recall that in the Newton--Raphson method, we start with
an initial approximation $x_0$ and generate a sequence of
approximations $\{x_n\}^\infty_{n=0}$ to $\alpha$ through
 \begin{equation}\label{eq6}
 x_{n+1}=x_n-\frac{f(x_n)}{f'(x_n)},\quad n=0,1,\ldots\ .
 \end{equation}
We also recall that, when $f(x)$ is twice continuously
differentiable in a closed  interval $I$ that includes $\alpha$,
and $f'(\alpha)\neq 0$, this method has order $2$. As such, the
Newton--Raphson method is extremely effective . To avoid computing
$f'(x)$ [note that $f'(x)$ may not always be available or may be
costly to compute], and to preserve the excellent convergence
properties of the Newton--Raphson method, we replace $f'(x_n)$ in
\eqref{eq6} by the approximation
$f[x_n,x_{n-1}]=[f(x_n)-f(x_{n-1})]/(x_n-x_{n-1})$. This results
in \eqref{eq4}, that is, in the secant method. The justification
for this approach is as follows:  When convergence takes place,
that is, when $\lim_{n\to\infty}x_n=\alpha$, the difference
$x_n-x_{n-1}$ tends to zero, and this implies that, as $n$
increases,  the accuracy of $f[x_n,x_{n-1}]$ as an approximation
to  $f'(x_n)$ increases as well.

In Section~\ref{se2} of this note, we consider in detail a
generalization of the second of the two approaches described above
using polynomial interpolation of degree $k$ with $k>1$. The
$(k+1)$-point iterative method that results from this
generalization turns out to be very effective. It is of order
higher than that of the secant method and requires only one
function evaluation per iteration. In Section~\ref{se3}, we
analyze this method and determine its order as well. In Section
\ref{se4}, we confirm our theory via a numerical example.

This paper is a slightly extended version of the paper \cite{Sidi:2008:GSM}.
The original version concerns the case in which $f^{(k+1)}(\alpha)\neq0$, while the extension concerns the special case in which $f^{(k+1)}(\alpha)=0$, which always occurs when $f(x)$ is a polynomial of degree at most $k$. In addition, we include a brief
discussion of the efficiency index for our method as  Section \ref{se5}.

\section{Generalization of secant method}\label{se2}
We start by discussing a known generalization of the secant method
(see, for example, Traub \cite[Chapters 4, 6, and
10]{Traub:1964:IMS}). In this generalization, we approximate
$f(x)$ by the polynomial of interpolation $p_{n,k}(x)$, where
$p_{n,k}(x_i)=f(x_i),$ $i=n,n-1,\ldots,n-k,$ assuming that
$x_0,x_1,\ldots,x_n$  have all been computed. Following that, we
determine $x_{n+1}$ as a zero of $p_{n,k}(x)$, provided a real
solution to $p_{n,k}(x)=0$ exists. Thus, $x_{n+1}$ is the solution
to a polynomial equation of degree $k$. For $k=1$, what we have is
nothing but the secant method. For $k=2$, $x_{n+1}$ is one of the
solutions to a quadratic equation, and the resulting method is
known as the method of M\"{u}ller. Clearly, for $k\geq 3$, the
determination of $x_{n+1}$ is not easy.

This difficulty prompts us to consider the second approach to the
secant method we discussed in Section \ref{se1}, in which we
replaced $f'(x_n)$ by the slope of the straight line through the
points $(x_n,f(x_n))$ and $(x_{n-1},f(x_{n-1}))$, that is, by the
derivative (at $x_n$) of the (linear) interpolant to $f(x)$ at
$x_n$ and $x_{n-1}$. We generalize this approach by replacing
$f'(x_n)$ by $p'_{n,k}(x_n)$,  the derivative at $x_n$ of the
 polynomial $p_{n,k}(x)$ interpolating $f(x)$ at the points $x_{n-i}$,
 $i=0,1,\ldots,k$, mentioned in the preceding paragraph, with $k\geq 2$.
Because $p_{n,k}(x)$ is a better approximation to $f(x)$ in the
neighborhood of $x_n$, $p'_{n,k}(x_n)$ is a better approximation
to $f'(x_n)$ when $k\geq 2$ than $f[x_n,x_{n-1}]$ used in the
secant method. In addition, just as the secant method, the new
method computes the function $f(x)$ only once per iteration step,
the computation being that of $f(x_n)$.
 Thus, the new method is described by the following $(k+1)$-point
 iteration:

\begin{equation}\label{eq:111}
x_{n+1}=x_n-\frac{f(x_n)}{p'_{n,k}(x_n)}, \quad n=k,k+1,\ldots,
\end{equation}
with $x_0,x_1,\ldots,x_k$ as initial approximations to be provided
by the user. Of course, with $k$ fixed, we can start with $x_0$
and $x_1$,  compute $x_2$ via the method we have described (with
$k=1$, namely via the secant method), compute $x_3$ via the method
we have described (with $k=2$), and so on, until we have completed
the list $x_0,x_1,\ldots,x_k$.

We now turn to the computational aspects of this method. What we
need is a fast method for computing $p'_{n,k}(x_n)$. For this, we
write $p_{n,k}(x)$ in Newtonian form as follows:
\begin{equation}\label{eq:151}
 p_{n,k}(x)=f(x_n)+\sum^{k}_{i=1}f[x_n,x_{n-1},\ldots,x_{n-i}]
\prod^{i-1}_{j=0}(x-x_{n-j}).
\end{equation}
Here $f[x_i,x_{i+1},\ldots,x_m]$ are divided differences of
$f(x)$, and  we recall that they can be defined recursively via
\begin{equation} \label{dda}
 f[x_i]=f(x_i);\quad f[x_i,x_j]=\frac{f[x_i]-f[x_j]}{x_i-x_j}, \quad x_i\neq x_j,
 \end{equation}
and, for $m>i+1$, via
\begin{equation} \label{ddb}
f[x_i,x_{i+1},\ldots,x_m]=\frac{f[x_i,x_{i+1},\ldots,x_{m-1}]-f[x_{i+1},x_{i+2},\ldots,x_m]}
{x_i-x_m},\quad x_i\neq x_m.
\end{equation}
We also recall that $f[x_i,x_{i+1},\ldots,x_m]$ is a symmetric
function of its arguments, that is, it has the same value under
any permutation of $\{x_i,x_{i+1},\ldots,x_m\}$. Thus, in
\eqref{eq:151},
 $$ f[x_n,x_{n-1},\ldots,x_{n-i}]=f[x_{n-i},x_{n-i+1},\ldots,x_n].$$
 In addition, when  $f\in \mathbb{C}^m(I)$, where $I$ is an open interval containing
 the points $z_0,z_1,\ldots,z_m$, whether these are distinct or not, there holds
 $$ f[z_0,z_1,\ldots,z_m]=\frac{f^{(m)}(\xi)}{m!}\quad \text{for
 some $\xi\in (\min\{z_i\},\max\{z_i\})$.}$$

 \sloppypar Going back to  \eqref{eq:151},  we note that  $p_{n,k}(x)$ there is computed
 by ordering    the
$x_i$ as $x_n,x_{n-1},\ldots,x_{n-k}$. This ordering enables us to
compute  $p'_{n,k}(x_n)$ easily. Indeed, differentiating
$p_{n,k}(x)$ in  \eqref{eq:151}, and letting $x=x_n$, we obtain

\begin{equation}\label{eq:152}
p'_{n,k}(x_n)=f[x_n,x_{n-1}]+\sum^{k}_{i=2}f[x_n,x_{n-1},\ldots,x_{n-i}]
\prod^{i-1}_{j=1}(x_n-x_{n-j}).
\end{equation}

In addition, note that the relevant divided difference table need
not be computed anew each iteration; what is needed is adding a
new diagonal (from  the south-west to north-east) at the bottom of
the existing table. To make this point clear, let us look at the
following example: Suppose $k=3$ and we have computed $x_i$,
$i=0,1,\ldots,7.$ To compute $x_8$, we use the divided difference
table in Table~\ref{table1}. Letting $f_{i,i+1,\ldots,m}$ stand
for $f[x_i,x_{i+1},\ldots,x_m]$, we have
\begin{table}[t]
\begin{center}
$$
\begin{matrix}
x_0&f_0&&&\\
&&f_{01}&&\\
x_1&f_1&&f_{012}&\\
&&f_{12}&&f_{0123}\\
x_2&f_2&&f_{123}&\\
&&f_{23}&&f_{1234}\\
x_3&f_3&&f_{234}&\\
&&f_{34}&&f_{2345}\\
x_4&f_4&&f_{345}&\\
&&f_{45}&&f_{3456}\\
x_5&f_5&&f_{456}&\\
&&f_{56}&&f_{4567}\\
x_6&f_6&&f_{567}&\\
&&f_{67}&&\\
x_7&f_7&&&\\

\end{matrix}
$$
\caption{\label{table1} \small Table of divided differences over
$\{x_0,x_1,\ldots,x_7\}$ for use to compute $x_8$ via $p_{7,3}(x)$
in \eqref{eq:111}. Note that $f_{i,i+1,\ldots,m}$ stands for
$f[x_i,x_{i+1},\ldots,x_m]$ throughout.}
\end{center}
\end{table}

$$ x_8=x_7-\frac{f(x_7)}{p_{7,3}(x_7)}=
x_7-\frac{f_7}{f_{67}+f_{567}(x_7-x_6)+f_{4567}(x_7-x_6)(x_7-x_5)}.$$
 To compute $x_9$, we will need the divided differences
$f_8,f_{78},f_{678},f_{5678}$. Computing first $f_8=f(x_8)$ with
the newly computed $x_8$, the rest of these divided differences
can be computed from the bottom diagonal of Table~\ref{table1} via
the recursion relations
$$f_{78}=\frac{f_7-f_8}{x_7-x_8},\quad f_{678}=\frac{f_{67}-f_{78}}{x_6-x_8},
\quad f_{5678}=\frac{f_{567}-f_{678}}{x_5-x_8},\quad\text{in this
order,}$$
 and appended to the bottom of Table~\ref{table1}.
Actually, we can do even better: Since we need only the bottom
diagonal of Table \ref{table1} to compute $x_8$, we need to save
only this diagonal, namely, only the entries
$f_7,f_{67},f_{567},f_{4567}.$ Once we have computed $x_8$ and
$f_8=f(x_8)$, we can overwrite $f_7,f_{67},f_{567},f_{4567}$ with
$f_8,f_{78},f_{678},f_{5678}$. Thus, in general, to be able to
compute $x_{n+1}$ via \eqref{eq:111}, after $x_n$ has been
determined, we need to store only the entries
$f_n,f_{n-1,n},\ldots,f_{n-k,n-k+1,\ldots,n-1,n}$ along with
$x_n,x_{n-1},\ldots,x_{n-k}$.

\section{Convergence analysis}\label{se3} We now turn to
the analysis of the sequence $\{x_n\}^{\infty}_{n=0}$ that is
generated via \eqref{eq:111}. Since we already know everything
concerning the case $k=1$, namely, the secant method, we treat the
case $k\geq 2$. The following theorem gives the main convergence
result for the generalized secant method.

\begin{theorem}\label{th1}
Let $\alpha$ be the solution to the equation $f(x)=0$. Assume
$f\in C^{k+1}(I)$, where $I$  is an open interval containing
$\alpha$, and assume also that  $f'(\alpha)\neq 0$, in addition to
$f(\alpha)=0$. Let $x_0,x_1,\ldots,x_k$ be distinct initial
approximations to $\alpha$, and generate $x_n$,
$n=k+1,k+2,\ldots,$ via
 \begin{equation} \label{eq:gs}
x_{n+1}=x_n-\frac{f(x_n)}{p'_{n,k}(x_n)},\quad
n=k,k+1,\ldots,\end{equation} where $p_{n,k}(x)$ is the polynomial
of interpolation to $f(x)$ at the points
$x_n,x_{n-1},\ldots,x_{n-k}$.
Then, provided $x_0,x_1,\ldots,x_k$ are in $I$ and sufficiently close to
$\alpha$, we have the following  cases:
\begin{enumerate}
\item
 If $f^{(k+1)}(\alpha)\neq0$, the sequence $\{x_n\}$ converges to $\alpha$, and

\begin{equation}\label{eq:130}
\lim_{n\to\infty}\frac{\epsilon_{n+1}}{\prod^k_{i=0}\epsilon_{n-i}}=
\frac{(-1)^{k+1}}{(k+1)!}\,\frac{f^{(k+1)}(\alpha)}{f'(\alpha)}\equiv
L;\quad \epsilon_n=x_n-\alpha\quad \forall n. \end{equation}
 The order of convergence is $s_k$, $1<s_k<2$, where $s_k$ is the only
 positive root of the equation $s^{k+1}=\sum_{i=0}^k s^i$ and
 satisfies
\begin{equation}\label{eq:132} 2-2^{-k-1}e< s_k<2-2^{-k-1}\quad\text{for $k\ge 2$;} \quad s_k<s_{k+1};\quad
\lim_{k\to\infty}s_k=2,\end{equation} where  $e$ is the base of
natural logarithms, and
\begin{equation}\label{eq:135}
\lim_{n\to\infty}\frac{|\epsilon_{n+1}|}{~~~|\epsilon_{n}|^{s_k}}=|L|^{(s_k-1)/k}.\end{equation}
\item
 If $f(x)$ is a polynomial of degree at most $k$, the sequence $\{x_n\}$ converges to $\alpha$, and
\begin{equation}\label{eq:177}
\lim_{n\to\infty}\frac{\epsilon_{n+1}}{\epsilon_n^2}=
\frac{f''(\alpha)}{2f'(\alpha)};\quad \epsilon_n=x_n-\alpha\quad \forall n. \end{equation}
Thus $\{x_n\}$ converges  of order $2$ if $f''(\alpha)\neq0$, and of  order  greater than $2$ if $f''(\alpha)=0$.
\end{enumerate}
\end{theorem}

\noindent{\bf Remark.} Note that, in part 1 of  Theorem \ref{th1},
$$s_1\dot{=}1.618, \ \ s_2\dot{=}1.839,\ \ s_3\dot{=}1.928, \ \ s_4\dot{=}1.966, \ \
s_5\dot{=}1.984, \ \ s_6\dot{=}1.992, \ \ s_7\dot{=}1.996,\ \ \text{etc.} $$
rounded to four significant figures. (Recall that $s_1$ is the order of the secant method.)\\

\noindent{\bf Proof.} Below, we shall  use the short-hand notation
$$
\text{int}(a_1,\ldots,a_m)=(\min\{a_1,\ldots,a_m\},\max\{a_1,\ldots,a_m\}).$$

We start by deriving a closed-form expression for the error in
$x_{n+1}$.  Subtracting $\alpha$ from both sides of \eqref{eq:gs},
and noting that
$$ f(x_n)=f(x_n)-f(\alpha)=f[x_n,\alpha](x_n-\alpha),$$
we have
\begin{equation}\label{eq:112}
x_{n+1}-\alpha=\bigg(1-\frac{f[x_n,\alpha]}{p'_{n,k}(x_n)}\bigg)(x_n-\alpha)=
\frac{p'_{n,k}(x_n)-f[x_n,\alpha]}{p'_{n,k}(x_n)}\,(x_n-\alpha).
\end{equation}
 We now note that
$$
p'_{n,k}(x_n)-f[x_n,\alpha]=\big[p'_{n,k}(x_n)-f'(x_n)\big]+
\big[f'(x_n)-f[x_n,\alpha]\big],$$   which, by
\begin{align*}
f'(x_n)-f[x_n,\alpha]&=f[x_n,x_n]-f[x_n,\alpha]\\
&=f[x_n,x_n,\alpha](x_n-\alpha)\\&=
\frac{f^{(2)}(\eta_n)}{2!}(x_n-\alpha)\quad \text{for some
$\eta_n\in\text{int}(x_n,\alpha)$,}\end{align*} and
\begin{align} \label{eq:107}
f'(x_n)-p'_{n,k}(x_n)&=f[x_n,x_n,x_{n-1},\ldots,x_{n-k}]\prod^k_{i=1}(x_n-x_{n-i})\nonumber \\
&=\frac{f^{(k+1)}(\xi_{n})}{(k+1)!}\prod^k_{i=1}(x_n-x_{n-i})
\quad \text{for some
$\xi_{n}\in\text{int}(x_n,x_{n-1},\ldots,x_{n-k})$,}
\end{align}
 becomes
\begin{equation}\label{eq:109}p'_{n,k}(x_n)-f[x_n,\alpha]=-\frac{f^{(k+1)}(\xi_{n})}{(k+1)!}
\prod^k_{i=1}(\epsilon_n-\epsilon_{n-i})+\frac{f^{(2)}(\eta_n)}{2!}\epsilon_n.\end{equation}
Substituting  \eqref{eq:107} and \eqref{eq:109}   in
\eqref{eq:112}, and letting

\begin{equation}\label{eq:128}
\widehat{D}_n=-\frac{f^{(k+1)}(\xi_{n})}{(k+1)!}\quad\text{and}\quad
\widehat{E}_n=\frac{f^{(2)}(\eta_n)}{2!},\end{equation} we finally
obtain

\begin{equation}\label{eq:123}
\epsilon_{n+1}=C_n\epsilon_n;\quad
C_n\equiv\frac{p'_{n,k}(x_n)-f[x_n,\alpha]}{p'_{n,k}(x_n)}=\frac{\widehat{D}_n\prod^k_{i=1}(\epsilon_n-\epsilon_{n-i})
+\widehat{E}_n
\epsilon_n}{f'(x_n)+\widehat{D}_n\prod^k_{i=1}(\epsilon_n-\epsilon_{n-i})}.\end{equation}

We now prove that convergence takes place.  Let $M_s=\max_{x\in
I}|f^{(s)}(x)|/s!,$ $s=1,2,\ldots,$ and choose the interval
$I=(\alpha-t/2,\alpha+t/2)$ sufficiently small to ensure that
$m_1=\min_{x\in I}|f'(x)|>0$ and $m_1>2M_{k+1}t^k+M_2t/2$. This is
possible since $\alpha\in I$ and $f'(\alpha)\neq 0.$ It can now be
shown that, provided $x_{n-i}$, $i=0,1,\ldots,k,$ are all in $I$,
there holds

\begin{align*}|C_n|&\leq
\frac{M_{k+1}\prod^k_{i=1}|\epsilon_n-\epsilon_{n-i}|+M_2|\epsilon_n|}
{m_1-M_{k+1}\prod^k_{i=1}|\epsilon_n-\epsilon_{n-i}|}  \\
&\leq
\frac{M_{k+1}\prod^k_{i=1}(|\epsilon_n|+|\epsilon_{n-i})|+M_2|\epsilon_n|}
{m_1-M_{k+1}\prod^k_{i=1}(|\epsilon_n|+|\epsilon_{n-i}|)}\leq
\overline{C},\end{align*} where
$$
\overline{C}\equiv\frac{M_{k+1}t^k+M_2t/2}{m_1-M_{k+1}t^k}<1.$$
Consequently, by \eqref{eq:123}, $|\epsilon_{n+1}|<|\epsilon_n|$,
which implies that $x_{n+1}\in I$, just like $x_{n-i}$,
$i=0,1,\ldots,k.$ Therefore, if  $x_0,x_1,\ldots,x_k$ are chosen
in $I$, then $|C_n|\leq \overline{C}<1$ for all $n\geq k$, hence
$\{x_n\}\subset I$  and $\lim_{n\to\infty}x_n=\alpha$.

As for \eqref{eq:130} when $f^{(k+1)}(\alpha)\neq0$, we proceed as follows: By the fact that
$\lim_{n\to\infty}x_n=\alpha$,  we first note that

$$ \lim_{n\to\infty}p'_{n,k}(x_n)=f'(\alpha)=\lim_{n\to\infty}f[x_n,\alpha],$$
and thus $\lim_{n\to\infty}C_n=0$. This means that
$\lim_{n\to\infty}(\epsilon_{n+1}/\epsilon_n)=0$ and,
equivalently, that $\{x_n\}$ converges of order greater than 1.
  As a result,

$$\lim_{n\to\infty}(\epsilon_n/\epsilon_{n-i})=0\quad\text{for all $i\geq
1$,}$$ and

$$  \epsilon_n/\epsilon_{n-i}=o(\epsilon_n/\epsilon_{n-j})\quad\text{as
$n\to\infty$,\quad for $j<i$.}$$
 Consequently, expanding in \eqref{eq:123} the product $\prod^k_{i=1}(\epsilon_n-\epsilon_{n-i})$, we have

\begin{align}\label{eq:200}\prod^k_{i=1}(\epsilon_n-\epsilon_{n-i})&=\prod^k_{i=1}\bigg(-\epsilon_{n-i}
[1-\epsilon_n/\epsilon_{n-i}]\bigg)
 \nonumber \\
&= (-1)^k
\bigg(\prod^k_{i=1}\epsilon_{n-i}\bigg)[1+O(\epsilon_n/\epsilon_{n-1})]\quad\text{as
$n\to\infty$.}\end{align}
 Substituting \eqref{eq:200}
in \eqref{eq:123},  and defining

\begin{equation}\label{eq:271}
D_n=\frac{\widehat{D}_n}{p'_{n,k}(x_n)}, \quad
E_n=\frac{\widehat{E}_n}{p'_{n,k}(x_n)},\end{equation}  we obtain

\begin{equation}\label{eq:231} \epsilon_{n+1}=(-1)^k
D_n\bigg(\prod^k_{i=0}\epsilon_{n-i}\bigg)[1+O(\epsilon_n/\epsilon_{n-1})]+E_n\epsilon_n^2\quad
\text{as $n\to\infty$.}
\end{equation}
Dividing both sides of \eqref{eq:231} by
$\prod^k_{i=0}\epsilon_{n-i}$, and defining

\begin{equation}\label{eq:233}
\sigma_n=\frac{\epsilon_{n+1}}{\prod^k_{i=0}\epsilon_{n-i}},
\end{equation}
 we have

\begin{equation}\label{eq:235}
 \sigma_n=(-1)^kD_n[1+O(\epsilon_n/\epsilon_{n-1})] +E_n \sigma_{n-1}\epsilon_{n-k-1}\quad \text{as
$n\to\infty$.}
\end{equation}
Now, \begin{equation}\label{eq:237}
\lim_{n\to\infty}D_n=-\frac{1}{(k+1)!}\frac{f^{(k+1)}(\alpha)}{f'(\alpha)},\quad
\lim_{n\to\infty}E_n=\frac{f^{(2)}(\alpha)}{2f'(\alpha)}.\end{equation}
Because  $\lim_{n\to\infty}D_n$ and $\lim_{n\to\infty}E_n$ are
finite, $\lim_{n\to\infty}(\epsilon_n/\epsilon_{n-1})=0$, and
$\lim_{n\to\infty}\epsilon_{n-k-1}=0$,  it follows that there
exist a positive integer $N$ and positive constants $\beta<1$ and
$D$,  with $|E_n\epsilon_{n-k-1}|\leq \beta$ when $n\geq N$, for
which \eqref{eq:235} gives

\begin{equation}\label{eq:241}
 |\sigma_n|\leq D+\beta|\sigma_{n-1}|\quad\text{for all
$n\geq N$.}\end{equation} Using \eqref{eq:241}, it is easy to show
that
$$|\sigma_{N+s}|\leq
D\frac{1-\beta^s}{1-\beta}+\beta^s|\sigma_N|,\quad s=1,2,\ldots,$$
which, by the fact that $\beta<1$, implies that $\{\sigma_n\}$ is
a bounded sequence. Making use of this fact, we
 have $\lim_{n\to\infty}E_n \sigma_{n-1}\epsilon_{n-k-1}=0$. Substituting this
 in \eqref{eq:235}, and  invoking \eqref{eq:237}, we next obtain
$\lim_{n\to\infty}\sigma_n=(-1)^k\lim_{n\to\infty}D_n=L$, which is
precisely \eqref{eq:130}.

That the order of the method is $s_k$, as defined in the statement
of the theorem, follows from
 \cite[Chapter 3]{Traub:1964:IMS}.  A  weaker version can be proved by
letting $\sigma_n=L$ for all $n$ and showing that
$|\epsilon_{n+1}|=Q|\epsilon_n|^{s_k}$ is possible for $s_k$ a
solution to the equation $s^{k+1}=\sum^k_{i=0}s^i$ and
$Q=|L|^{(s_k-1)/k}$. The  proof of this is easy and is left  to
the reader. This completes the proof of part 1 of the theorem.

When $f(x)$ is a polynomial of degree at most $k$, we first observe that $f^{(k+1)}(x)=0$ for all $x$, which implies that $p_{n,k}(x)= f(x)$ for all $x$, hence also
$p'_{n,k}(x)= f'(x)$ for all $x$.  Therefore, we have that
$p'_{n,k}(x_n)=f'(x_n)$ in the recursion  of \eqref{eq:gs}. Consequently,  \eqref{eq:gs} becomes
$$ x_{n+1}=x_n-\frac{f(x_n)}{f'(x_n)}, \quad n=k,k+1,\ldots, $$
which is the recursion for the Newton--Raphson method.  Thus,
\eqref{eq:177} follows. This completes the  proof of part 2 of the theorem.
\hfill $\blacksquare$

\section{A numerical example}\label{se4}

 We apply the method described in Sections \ref{se2} and
\ref{se3} to the solution of the equation $f(x)=0$, where
$f(x)=x^3-8$, whose solution is $\alpha=2.$ We take $k=2$ in our
method.   We also chose $x_0=0$ and $x_1=6$, and compute $x_2$ via
one step of the secant method, namely,
\begin{equation} \label{eq:5-1}
x_2=x_1-\frac{f(x_1)}{f[x_0,x_1]}.\end{equation}
 Following that,
we compute $x_3,x_4,\ldots,$ via
\begin{equation} \label{eq:5-2}
 x_{n+1}=x_n-\frac{f(x_n)}
{f[x_n,x_{n-1}]+f[x_n,x_{n-1},x_{n-2}](x_n-x_{n-1})},\quad
n=2,3,\ldots\ .\end{equation}

 Our computations were done in
quadruple-precision arithmetic (approximately 35-decimal-digit
accuracy), and they are given in Table \ref{table10}. Note that in
order to verify the theoretical results concerning iterative
methods of order greater that unity, we need to use computer
arithmetic of high precision (preferably, of variable precision,
if available) because the number of correct significant decimal
digits increases dramatically from one iteration to the next as we
are approaching the solution.

\begin{table}[tb]
\begin{center}
$$
\begin{array}{||c||c|c|c|c||}
\hline
 n& x_n &\epsilon_n&\displaystyle
 \frac{\epsilon_{n+1}}{\epsilon_n\epsilon_{n-1}\epsilon_{n-2}}
 & \displaystyle\frac{\log|\epsilon_{n+1}/\epsilon_n|}{\log|\epsilon_{n}/\epsilon_{n-1}|}\\
\hline\hline
   0&    5.00000000000000000000000000000000000D+00&      ~~3.000D+00&           &\\
  1&    4.00000000000000000000000000000000000D+00&      ~~2.000D+00&      &      1.515\\
  2&    3.08196721311475409836065573770491792D+00&      ~~1.082D+00&      ~~0.0441&      2.164\\
  3&    2.28621882971781130732266803773062580D+00&      ~~2.862D-01&      ~~0.1670&      2.497\\
  4&    2.01034420943787831264152973172014271D+00&      ~~1.034D-02&     -0.6370&      1.182\\
  5&    1.99979593345266992578358353656798415D+00&     -2.041D-04&     -0.1196&      2.024\\
  6&    2.00000007223139333059960671366229837D+00&     ~~7.223D-08&     -0.1005&      1.934\\
  7&    2.00000000000001531923884491258853168D+00&      ~~1.532D-14&     -0.0838&      1.784\\
  8&    2.00000000000000000000000001893448134D+00&      ~~1.893D-26&     * &*            \\
  9&    2.00000000000000000000000000000000000D+00&      ~~0.000D+00&  *  &*     \\
  \hline

\end{array}
$$

\caption{\label{table10} Results obtained by applying the
generalized secant method with $k=2$, as shown in \eqref{eq:5-1}
and \eqref{eq:5-2}, to the equation $x^3-8=0$.}

\end{center}
\end{table}

From Theorem \ref{th1},

$$
\lim_{n\to\infty}\frac{\epsilon_{n+1}}{\epsilon_{n}\epsilon_{n-1}\epsilon_{n-2}}=
\frac{(-1)^{3}}{3!}\,\frac{f^{(3)}(2)}{f'(2)}=-\frac{1}{12}=-0.08333\cdots$$
and $$ \lim_{n\to\infty}\frac{\log|\epsilon_{n+1}/\epsilon_{n}|}
{\log|\epsilon_{n}/\epsilon_{n-1}|}=s_2=1.83928\cdots,$$ and these
seem to be confirmed in Table \ref{table10}. Also, $x_9$
 should have a little under  50 correct significant
figures, even though we do not see this in  Table \ref{table10}
due to the fact that the arithmetic we have used to generate Table
\ref{table10} can provide an accuracy of at most 35 digits
approximately.

\section{Discussion of efficiency index of the method} \label{se5}
 We recall that, for methods that converge superlinearly, that is, with order strictly greater than $1$, a good measure of their effectiveness  is the so-called {\em efficiency index}, a concept introduced  originally by Ostrowski \cite{Ostrowski:1960:SES}. (See Traub \cite[pp. 11--13, 260--264]{Traub:1964:IMS} for more on this subject.)   If  an iterative method for solving $f(x)=0$  that requires $p$ evaluations of $f(x)$ (and its derivatives, assuming that their cost is about  the same), has  order $s>1$, the efficiency index $EI$ of the method is defined as $EI=s^{1/p}$. Figuratively speaking, $EI$ measures the order of the method {\em per function evaluation}.
 Thus, we may conclude that, the larger $EI$, the more effective the iterative method, {\em irrespective of its order}. In comparing methods, we should examine their performance after we  have done  a fixed  number of function evaluations, this number being {\em the same} for all methods. In other words, it makes sense to compare methods that have {\em equal} costs. The details of this line of thought follow:

  Consider    two iterative methods M1 and M2 applied to the equation $f(x)=0$, and let $m_1$ and $m_2$ be the number of function evaluations per iteration for M1 and M2, respectively.  Starting with $x^{(1)}_0=x^{(2)}_0$,  let the sequences of approximations $\{x^{(1)}_n\}^\infty_{n=0}$ and $\{x^{(2)}_n\}^\infty_{n=0}$ be  generated by M1 and M2, respectively.
  Then, for each integer $q=1,2,\ldots,$  we should compare the approximations $x^{(1)}_{qm_2}$ and  $x^{(2)}_{qm_1}$.
  Note that, the computation of  $x^{(1)}_{qm_2}$ starting from $x^{(1)}_{(q-1)m_2}$
  entails $m_1m_2$ function evaluations and so does
   the computation of  $x^{(2)}_{qm_1}$ starting from $x^{(2)}_{(q-1)m_1}$.

 In a fundamental paper by Kung and Traub \cite{Kung:1974:OOO}, it is conjectured that the order of a multipoint iterative method   without memory that uses $p$ function evaluations  may not exceed $2^{p-1}$. This paper contains two such families that use $p$ function evaluations and are of order $2^{p-1}$.
  Wo\'{z}niakowski \cite{Wozniakowski:1976:MOM} has proved for some classes of
 multipoint iterative methods   without memory that the order $2^{p-1}$ cannot be exceeded without more information. From this, it is clear that the efficiency index of such methods is at most $2^{(p-1)/p}=2^{1-1/p}<2$. In view of this discussion, we make a few comments on the  efficiency index of our method next.

   The efficiency index of the generalized secant method developed in this paper is $EI_k=s_k$ for each $k=1,2,\ldots,$  because $p=1$ for every $k$. In addition, because $\lim_{k\to\infty}s_k=2$, we have $\lim_{k\to\infty}IE_k=2$ as well.
  Actually, even with  very small $k$, we are able to come quite close to this limit.   For example, $s_7=1.9960\cdots$ and $s_{10}=1.9995\cdots$.

  Over the years, many sophisticated iterative methods  with and without memory that do not use derivatives of $f(x)$ and that have  high orders  have been developed. It is not our purpose here to review these methods; we refer the reader to the papers by   D\v{z}uni\'{c} \cite{Dzunic:2013:ETP}, D\v{z}uni\'{c} and Petkovi\'{c}
  \cite{Dzunic:2012:GMR}, \cite{Dzunic:2014:GBM},
  Chun and Neta \cite{Chun:2016:CSF},\cite{Chun:2019:CSM},
    and to the book by Petkovi\'{c} et al. \cite{Petkovic:2013:MMS}, for example,  and to the bibliographies of these
  publications.
     We only would like to comment that the many methods that we have studied have efficiency indices that are strictly less than $2$ despite their  high order.
  This may suggest that the method of this paper may be as  useful  a tool for solving nonlinear equations with simple zeros as other methods that have  orders
  much higher than $2$.


\end{document}